 \documentclass[11pt]{article}
\usepackage{amssymb, amsthm, amsmath, amscd}
\newtheorem{thm}{Theorem}[section]
\newtheorem{lem}[thm]{Lemma}
\newtheorem{prop}[thm]{Proposition}
\newtheorem{cor}[thm]{Corollary}

\theoremstyle{definition}\newtheorem{df}[thm]{Definition}
\theoremstyle{definition}\newtheorem{rem}[thm]{Remark}
\theoremstyle{definition}

\renewcommand{\phi}{\varphi}

\newcommand{\Q}{\mathbb{Q}}

\newcommand{\C}{\mathbb{C}}

\newcommand{\hm}{homomorphism}
\newcommand{\dt}{\delta}
\newcommand{\ep}{\epsilon}
\newcommand{\la}{\langle}
\newcommand{\ra}{\rangle}
\newcommand{\andeqn}{\,\,\,{\rm and}\,\,\,}

\newcommand{\CA}{$C^*$-algebra}
\newcommand{\SCA}{$C^*$-subalgebra}

\newcommand{\af}{{\alpha}}

\newcommand{\beq}{\begin{eqnarray}}
\newcommand{\eneq}{\end{eqnarray}}
\newcommand{\tforal}{\,\,\,\text{for\,\,\,all}\,\,\,}
\newcommand{\tand}{\,\,\,\text{and}\,\,\,}

\title{Cuntz semigroups of  \CA s of stable rank one and projective Hilbert modules}
\author{Huaxin Lin
 }
\date{}

\begin{document}

\maketitle

\begin{abstract}
Let $A$ be a simple \CA\, of stable rank one and let $p$ and $q$ be two
$\sigma$-compact open projections. It is proved that there is a
continuous path of unitaries in ${\tilde A}$ which connects open
sub-projections of $p$ which is  compactly contained in $p$ to those in
$q.$  It is also shown that every Hilbert module is projective in
the category whose morphisms are bounded module maps with adjoints.
A discussion of projective Hilbert modules (whose morphisms are
bounded module maps) is also given.
\end{abstract}

Recently the Cuntz semigroups of \CA s have attracted some
previously unexpected attention. The Cuntz relations for positive
elements in a \CA\, was introduced by J. Cuntz (see \cite{Cu1}). The
Cuntz semigroups, briefly, are semigroups  of equivalence classes of
positive elements in a \CA s. This  relation is similar to the Murry
and Von Neumann equivalence relation for projections. The renew
interests in the Cuntz semigroups probably begins with Toms's
example (\cite{Toms} ) which shows that two  unital simple
AH-algebras with the same  traditional Elliott invariant may have
different Cuntz semigroups.  It is a hope of many that  the Cuntz semigroups may be
used in the classification of amenable \CA s. This note limits
itself
 to the clarification of  a couple issues related to the Cuntz
 semigroups and its relation with Hilbert modules.

While the Cuntz semigroups  may be  useful  tools to distinguish some
\CA s, they are  not necessarily easy to compute in general. One problem is that
the Cuntz semigroup is not a homotopy invariant.
Let $A$ be a \CA\, and let $f\in C([0,1], A)$ so that $f(t)\ge 0$
for all $t\in [0,1].$ One easily sees that $f(0)$ and $f(1)$ are
unlikely related  in the Cuntz relation. On the other hand, Cuntz
introduced several versions of the relation among positive elements in
\CA s. These relations also give equivalence relations among open
projections of \CA s. It will be presented, following a result of L.
G. Brown, if two $\sigma$-compact open projections are homotopy,
then they are actually equivalent in (a strong) Cuntz relation.

Another homotopy question is whether two positive elements are homotopy in a suitable sense
if they  are equivalent in
the sense of Cuntz.
Under the assumption that $A$ is simple and has stable rank one, it is shown in
this note that two $\sigma$-compact open projections are Cuntz
equivalent if and only if there is a continuous path of unitaries
$\{u(t): t\in [0,1)\}$ which connects these two open projections in
the sense that will be described in \ref{Thomp} and \ref{Rhomp}. In
particular, any pre-compact open subprojection
 (see \ref{compact}) of $p$ is unitarily equivalent to a pre-compact
open subprojection of $q.$

 Let $A$ be \CA\, and let $a, b\in A\otimes {\cal K}.$ Then
 $H_1=\overline{aA}$ and $H_2=\overline{bA}$ are two Hilbert $A$-modules. Suppose that
 $p_a$ and $p_b$ are range projections of $a$ and $b$ in $(A\otimes {\cal K})^{**}.$
 Then $p_a$ and $p_b$ are Cuntz equivalent (see \ref{Dcuntz}) if and only if
 $H_1$ and $H_2$ are isomorphic as Hilbert $A$-modules. So Hilbert
 modules and the Cuntz semigroups are closely related.  In this note, using 
 another result of L. G. Brown, 
 it is shown that $p_a$ is dominated by $p_b$ in the sense of Cuntz
 if and only if there is a bounded module map $T: H_2\to H_1$
 (which may not have an adjoint) whose range is dense in $H_1.$
 Projectivity
 of Hilbert modules have been recently brought into attention.  In
 the last section of this note, quite differently from the pure algebraic analogy,
 it is shown that every Hilbert module over a
 \CA\, $A$  is projective in the category of Hilbert  $A$-modules
 with bounded module maps with adjoints as morphisms. However, for
 Hilbert modules over \CA\, $A,$ sometime the category of Hilbert
 $A$-modules with bounded module maps (may or may not have adjoints)
 is also useful. To determine which Hilbert $A$-modules are
 projective in that category is more difficult. A discussion on
 this problem will  also be presented.


\section{The Cuntz Semigroups}

\begin{df}\label{Dcuntz}

Let $A$ be a \CA\, and let $a\in A_+.$ Denote by
$Her(a)=\overline{aAa}$ the hereditary \SCA\, of $A$ generated by
$a.$ Denote by $p_a$ the range projection of $a$ in $A^{**}.$ It is
an open projection of $A$ in $A^{**}.$ Dnote $Her(p)=pA^{**}p\cap
A=Her(a).$

Suppose that $a, b\in A_+.$ One writes $a\lessapprox b$ if there
exists $x\in A$ such that $x^*x=a$ and $xx^*\in Her(b).$ One writes
$a\lesssim b,$ if there exists a sequence $r_n\in A$ such that
$r_n^*br_n\to a$ in norm. If $a\lesssim b$ and $b\lesssim a,$ then
one writes $a\sim b.$ The relation ``$\sim $" is an equivalence
relation. The equivalence class represented by $a$ will be written
as $\langle a\rangle.$ Denote by $W(A)$ the equivalence classes of
positive elements in $M_{\infty}(A)$ with respect to ``$\sim $". So
$$
W(A)=\{\langle a\rangle : a\in M_{\infty}(A)\}.
$$
The semigroup $W(A)$ is called the Cuntz semigroup. 
One can also define the same relation in $A\otimes {\cal K}.$ The corresponding  semigroup 
is denoted by $Cu(A).$ 

  Let $p$ and $q$ be two open projections of $M_n(A)$ in $M_n(A)^{**}=(M_n(A^{**}))$ for some integer $n\ge 1$ 
  (or let $p$ and $q$ be two open projections of $A\otimes {\cal K},$ where ${\cal K}$ is the \CA\, of compact operators  on 
  $l^2$).  One says that $p$ and $q$ are Cuntz equivalent and writes $p\approx_{cu} q,$ if there exists a
  partial isometry $v \in M_n(A)^{**}$ (or $v\in(A\otimes {\cal K})^{**}$) such that
  $$
  v^*v=p,\,\, vv^*=q\andeqn vav^*\in Her(q)\tforal a\in Her(p).
  $$
  The relation ``$\approx_{cu}$" is also an equivalence relation.
  The equivalence class represented by $p$ will be denoted by $[p].$
  An open projection of $A$ is said to be $\sigma$-compact, if $p=p_a$ for
  some $a\in A_+.$
  Denote by $Co(A)$ the equivalence classes of $\sigma$-compact open projections of 
  $M_n(A)$ in $M_n(A)^{**}$ for all $n\ge 1.$ Denote by $Co(A\otimes {\cal K})$ the equivalence classes of 
  $\sigma$-compact open projections of 
  $A\otimes {\cal K}$ in $(A\otimes {\cal K})^{**}.$ 
   Note that $M_n(A\otimes {\cal K})\cong A\otimes {\cal K}.$ These also form semigroups.
  One write $[p]\le [q],$ if $p\approx_{cu} q'$ for some open
  projection $q'\le q.$

  Let $a, b\in M_{\infty}(A)$ be two positive elements. Then
  $[p_a]\le [p_b]$ if and only if $a\lessapprox b.$ One writes
  $p\sim_{cu} q$ if $[p]\le [q]$ and $[q]\le [p].$ This is also an
  equivalence relation. Denote by $\langle p\rangle $ the
  equivalence class represented by $p.$

  These relations were first introduced by Cuntz (see \cite{Cu1}) and the
  readers are referred to \cite{Cu1}, \cite{Cu2}, \cite{B1},  \cite{R1}  and \cite{KR} for more details.

\end{df}

There are significant differences between $W(A)$ and $Co(A)$  (and 
differences between ($Cu(A)$ and $Co(A\otimes {\cal K})$) in
general. An example that $W(A)\not=Co(A)$ for stably finite \CA\,
was given in \cite{BC}. Let $A$ be a purely infinite simple \CA\,
and let  $a,\, b\in A_+\setminus \{0\}.$ Then $\la a\ra =\la b \ra.$
Thus $Cu(A)$ contains only zero and one other element. It is not
quite useful in this case.  On the other hand, it follows from a
result of S. Zhang \cite{Zh} (see also Cor. 11 of  \cite{Lnopn})
that $Her(a)$ and $Her(b)$ are stable and isomorphic if neither
$p_a$ nor $p_b$ are in $A.$  In fact the isomorphism can be given by
an isometry. From this, one easily obtains the following.

\begin{prop}
Let $A$ be a purely infinite simple \CA\, Then
$$
Co(A)=V(A)\sqcup \{\infty\},
$$
where $V(A)$ is the Murry-Von Neumann equivalence classes of
projections in $M_{\infty}(A)$ and $\infty$ is represented by a
non-zero $\sigma$-compact  open projection which is not in $A.$
\end{prop}

However,  $W(A)$ and $Co(A)$ could be often the same.

\vspace{0.1in}

\begin{df}

Let $\ep>0.$ Define
$$
f_\ep(t)=\begin{cases} 0 & \text{if}\,\, t\in [0,\ep/2]\\
                       {\rm linear}\, & \text{if}\,\, t\in [\ep/2,
                       \ep]\\
                       1 & \text{if}\,\,\,t\in
                       [\ep,\infty).\end{cases}
                       $$

\end{df}

\begin{lem}\label{Ped}{\rm (G. K. Pedersen (Theorem 5 of \cite{P1}))}
Let $A$ be a \CA\, with stable rank one. Suppose that $x\in A.$
Then, for each $t\in (0, \|x\|],$ there is a unitary $u_t\in {\tilde
A}$ such that
$$
u_t^*p_tu_t=q_t,
$$
where $p_t$ is the open spectral projection of $|x|$ associated with
$(t, \|x\|]$ and $q_t$ is the open spectral projection of $|x^*|$
associated with $(t, \|x\|],$ respectively. Moreover,
$$
u_{t'}p_t=u_tp_t\andeqn u_{t'}p_tu_{t'}=q_t
$$
for all $0<t'<t<\|x\|.$
\end{lem}

\begin{proof}
Note, by Theorem 5 of \cite{P1}, since $A$ has stable rank one, for
each $t\in (0, \|x\|],$ there is a unitary $u_t\in {\tilde A}$ such
that $u_tp_t=vp_t,$ where $x=v|x|$ is the polar decomposition of $x$
in $A^{**}.$ Then
$$
u_{t'}p_t=u_{t'}p_{t'}p_t=vp_{t'}p_t=vp_t=u_tp_t
$$
for any $0<t'<t<\|x\|.$ In particular,
$$
u_{t'}p_tu_{t'}^*=u_tp_tu_t^*=q_t.
$$
\end{proof}

\begin{prop}\label{lesssim}
Let $A$ be a \CA\, with stable rank one and let $a, b\in A_+$ be two
positive element.
 Then the following are equivalent.

 {\rm (1)} $[p_a]\le [p_b],$

 {\rm (2)}  $a\lessapprox b,$

 {\rm (3)}  $a\lesssim b.$

\end{prop}

\begin{proof}
From the definition, (2) implies (3).  It is also known that (1) and
(2) are equivalent. It remains to show (3) implies (1). To simplify
notation, one may assume that $A$ is unital, $0\le a, b\le 1.$
 Suppose (3)
holds. Let $\{\ep_n\}$ be a strictly decreasing sequence of positive
numbers in $(0,1]$ such that $\sum_{n=1}^{\infty}\ep_n\le 1/2.$

By \cite{R1}, there is a unitary $w_1\in A$ such that
\beq\label{lesssim-1}
b_1=w_1f_{\ep_1/4}(a)w_1^*\le w_1f_{\ep_1/16}(a)w_1^*={\bar b_1}\in
Her(b).
\eneq
Note that ${\bar b_1}b_1=b_1.$ Let $x_1=w_1(f_{\ep_1/4}(a))^{1/2}.$
Then
\beq\label{lesssim-2}
x_1^*x_1=f_{\ep_1/4}(a)\andeqn x_1x_1^*=b_1.
\eneq
There is a unitary $w_2\in A$ such that
\beq\label{lessim-4}
w_2w_1f_{\ep_2/8}(a)w_1^* w_2^*=b_2\in Her(b).
\eneq
Denote $a_1=w_2w_1f_{\ep_1/4}(a)w_1w_2^*.$  Note that $a_1\in
Her(b)$ and $a_1b_2=a_1.$ Therefore
\beq\label{lessim-5}
(b_2-1)w_2w_1|x_1|w_1^*=((b_2-1)w_2w_1|x_1|w_1^*w_2^*)w_2=0.
\eneq
In other words,
\beq\label{lessim-6}
b_2w_2w_1|x_1|w_1^*=w_2w_1|x_1|w_1^*.
\eneq
Similarly,
\beq\label{lessim-7}
w_2w_1|x_1|w_1^*{\bar b_1}=w_2w_1|x_1|w_1^*.
\eneq
Therefore $y_1:=w_2w_1|x_1|w_1^*\in Her(b).$ Moreover,
\beq\label{lessim-8}
y_1^*y_1=w_1x_1^*x_1w_1^*=x_1x_1^*\andeqn
y_1y_1^*=w_2w_1f_{\ep_1/4}(a)w_1^*w_2^*.
\eneq

By applying \ref{Ped}, one obtains a unitary $z_1\in {\tilde
Her(b)}$ such that
\beq\label{lesssim-9}
z_1e_{1/4}(|y_1|)=w_2e_{1/4}(|y_1|)=w_2e_{1/4}(|x_1^*|),
\eneq
where $e_{1/4}(|y_1|)$ is the open spectral projection of
$|y_1|=|x_1^*|$ associated with $(1/4, 1].$ Note that,
\beq\label{lesssim-10}
e_{1/4}(|x_1^*|)&=& e_{1/4}(w_1f_{\ep_1/4}(a)w_1^*)\\
 &=&
w_1e_{1/4}(f_{\ep_1/4}(a))w_1^*\\\label{lesssim-10+}
&=&w_1e_{\dt_1}(a)w_1^*
\eneq
where $e_{1/4}(f_{\ep_1/4}(a))$ is the open spectral projection of
$f_{\ep_1/4}(a)$ associated with $(\ep_1/4, 1]$ and $e_{\dt_1}(a)$
is the open spectral projection of $a$ associated with $(\dt_1, 1]$
for some $\dt_1\in (\ep_1/4, 3\ep_1/8).$

By (\ref{lesssim-9}) and (\ref{lesssim-10+}),
\beq\label{lesssim-11-}
z_1^*w_2w_1e_{\dt_1}(a)&=&z_1^*(w_2w_1e_{\dt_1}(a)w_1^*)w_1\\
&=&z_1^*(z_1e_{1/4}(|x_1^*|))w_1=e_{1/4}(|x_1^*|)w_1\\\label{lesssim-11-+}
&=& w_1e_{\dt_1}(a)
\eneq

 Define $u_1=w_1$ and
$u_2=z_1^*w_2u_1=z_1^*w_2w_1,$ where one may view $z_1$ as a unitary
in $A.$
 It follows, for any $x\in {\overline{f_{\dt_1}(a)A}},$ by applying
 (\ref{lesssim-11-+}),
that
\beq\label{lesssim-11}
u_2x=u_2e_{\dt_1}(a)x=z_1^*w_2w_1e_{\dt_1}(a)=u_1e_{\dt_1}(a)x=u_1x
\eneq

Note also that $u_2yu_2^*\in Her(b)$ for all $y\in
Her(f_{\ep_2/8}(a)),$ and $f_{\ep_1}(a) \in Her(e_{\dt_1}(a)).$

By induction, for each $n,$ one obtains a sequence of unitaries
$u_n\in A$ such that
\beq\label{lesssim-12}
u_nyu_n^*&\in& Her(b)\tforal y\in Her(f_{\ep_n/8}(a))\andeqn\\
u_{n+1}x&=&u_nx \tforal x\in {\overline{f_{\ep_n}(a)A}}
\eneq
One then computes that
\beq\label{lesssim-13}
\lim_{n\to\infty}u_nx.
\eneq
converges for every $x\in {\overline{aA}},$ which defines  a unitary
isomorphism $U$ from ${\overline{aA}}$ into a Hilbert sub-module of
${\overline{bA}},$ which implies that $[p_a]\le [p_b].$

\end{proof}

\begin{rem}\label{EE}
There will be some discussion of Hilbert modules in the last section. A countably
generated Hilbert module may not have a countable dense set. Note
that in Proposition \ref{lesssim}, $A$ is not assumed to be
separable. The  argument above can also be used to prove the
following theorem which was proved in \cite{CEI}.
\end{rem}

\begin{thm}\label{TCEI}
Let $A$ be a \CA\, of stable rank one and let $a, b\in A_+.$ Then
the following are equivalent:

{\rm (1)} $[p_a]=[p_b];$

{\rm (2)} $\langle p_a\rangle =\langle p_b\rangle;$

{\rm (3)} $\langle a\rangle=\langle b\rangle.$

In particular, $Co(A)=W(A)$ and $Co(A\otimes {\cal K})=Cu(A).$ 
\end{thm}

\section{Homotopy}

It seems quite appropriate to begin with the following result of L.
G. Brown (\cite{Br1}).

\begin{prop}\label{Hp1}
Let $A$ be a \CA\, and let  $p$ and $q$ be two $\sigma$-compact open
projections of $A$ in $A^{**}.$  Suppose that there is a norm continuous path
$\{p(t): t\in [0,1]\}$ of $\sigma$-compact open projections such
that
\beq\label{Hp1-1}
p(0)=p\andeqn p(1)=q.
\eneq
Then $[p]=[q].$
\end{prop}

\begin{proof}
Let $0=t_0<t_1<t_2<\cdots t_n=1$ be a partition such that
$$
\|p(t_i)-p(t_{i-1})\|<1/2,\,\,\, i=1,2,...,n.
$$

It follows from (the proof of) 3.2 of \cite{Br1} that
\beq\label{Hp1-2}
[p(t_i)]=[p(t_{i-1})], \,\,\,i=1,2,....,n.
\eneq
Thus $[p]=[q].$

\end{proof}

\begin{df}\label{contained}
Let $A$ be a \CA. An open projection $q\in A^{**}$ is said to be pre-compact, if
there is a positive element $a\in A_+$ such that $qa=qa=q.$
If $p$ is another open projection and if there is $a\in Her(p)$ such that
$qa=aq=q,$ then one says that $q$  {\it is compactly contained}  in $p.$

\end{df}

\begin{lem}\label{cut}
Let $A$ be a \CA\, and let $a\in A_+.$ Suppose that $q\in A^{**}$ is
a projection for which $qa=aq=q.$ Then, $\chi_{(1, \|a\|]}(a)q=q.$
\end{lem}

\begin{proof}
Let $s>1.$ Denote by $p_{(s, \|a\|]}$ the spectral projection of $a$ in $A^{**}$ corresponding to the interval $(s, \|a\|].$
Since $qa=aq,$ $q$ commutes with $p_{(s, \|a\|]}.$ In particular,
$qp_{(s, \|a\|]}$ is a projection.  However, $qp_{(s, \|a\|]}=0.$ Otherwise
\beq\label{cut-1-1}
1<s\le \|qap_{(s, \|a\|]}\|\le \|aq\|=\|q\|=1.
\eneq

It follows that
\beq\label{cut-1-2}
qp_{(1, \|a\|]}=0.
\eneq

Let $0<r<1$ and $p_{[0,r]}=\chi_{[0, r]}(a)$ be the spectral projection corresponding
to the interval $[0,r].$ The assumption that $qa=aq=q$ implies that
\beq\label{cut-1}
p_{[0,r]}q=qp_{[0,r]}.
\eneq
It follows that
\beq\label{cut-2}
q=aq\le  rp_{[0,r]}q+ap_{(r, \|a\|]}q\le rp_{[0,r]}q+p_{(r, 1]}q\le p_{[0, r]}q+p_{(r, 1]}q=q.
\eneq
It follows that
\beq\label{cut-3}
rp_{[0,r]}q=p_{[0,r]}q.
\eneq
Therefore
\beq\label{cut-4}
p_{[0,r]}q=0.
\eneq
Since this holds for each $r\in (0,1),$ one concludes that
\beq\label{cut-5}
q=p_aq=\chi_{\{1\}}(a)q=q.
\eneq
\end{proof}

\begin{df}\label{sgmcomp}
Let $p$ be a $\sigma$-compact open projection of $A$ which is not in $A.$
Let $a\in Her(p)$ be a strictly positive element. Then $0$ must be a
limit point of ${\rm sp}(a).$ Let $t_n\in (0, \|a\|]$ be such that
$t_n\searrow 0.$ Let $p_{K,n}$ be the open spectral projection corresponding to  $(t_n, \|a\|].$ Then
$f_{t_n/2}(a)\ge p_{K,n},$ So $p_{K,n}$ is a sub-pre-compact open
projection of $p.$ Note that $\{p_{K,n}:n=1,2,...\}$ is increasing
and
$$
\lim_{n\to\infty}p_{K,n}=p
$$
in the strong operator topology in $A^{**}.$
Such a sequence $\{p_{K,n}\}$ is called a {\it pre-compact support} of
$p.$

\end{df}

In the proof of \ref{PK1}, \ref{Thomp} and \ref{compact}, the result of L. G. Brown
and G. K. Pedersen (3.6 of  \cite{BP}) that every hereditary \SCA\,
of a \CA\, of stable rank one has stable rank one will be used
without repeating this reference.

\begin{lem}\label{PK1}
Let $A$ be a simple \CA\, of stable rank one and let $x\in A.$
Suppose that $x=v|x|$ is the polar decomposition of $x$ in $A^{**}.$
Suppose also that $0$ is not an isolated point in ${\rm sp}(x).$
Then, for any $\dt>0,$ there is a unitary $u\in {\tilde A}$ with
$[u]=0$ in $K_1(A)$ such that
\beq\label{PK1-1}
up_t=vp_t\tforal t\in [\dt, \|x\|],
\eneq
where $p_t$ is the spectral projection of $|x|$ corresponding to
$(t, \|x\|].$
\end{lem}

\begin{proof}
It follows from  \ref{Ped} that there is a unitary $u_\dt\in {\tilde
A}$ such that
\beq\label{PK1-2}
u_\dt p_t=vp_t \tforal t\in [\dt, \|x\|].
\eneq
Since $0$ is not an isolated point in ${\rm sp}(x),$ there are
$0<t'<t''<\dt$ such that $p_{(t', t'')}\not=0,$ where $p_{(t',
t'')}$ is the spectral projection of $|x|$ corresponding to $(t', t'').$  Note that $p_{(t',
t'')}$ is an open projection of $A.$ Let $B=Her(p_{(t',t'')}).$ Then
$B$ has stable rank one. Since $A$ is also simple, the map
$K_1(B)\to K_1(A)$ induced by the inclusion is an isomorphism.
Therefore there is a unitary $v\in {\widetilde B}$ such that
$[v]=[u_\dt^*]$ in $K_1(A).$ One may write $v=z+\lambda,$ where
$z\in B$ and $\lambda\in \C.$ Let $\pi: {\widetilde{B}}\to \C$ be
the quotient map. Then $\pi(v)=\lambda.$ It follows that
$|\lambda|=1.$ Put
$v_1=\overline{\lambda}v={\overline{\lambda}}z+1.$  Note that
\beq\label{PK-3}
zp_t=0\tforal t\in (\dt, \|x\|]\andeqn [v_1]=[v]=[u_\dt^*]\,\,\,{\rm in}\,\,\, K_1(A).
\eneq
One may view $v_1$ as a unitary in ${\tilde A}.$   Now set $u=u_\dt
v_1.$ Then,
\beq\label{PK1-3}
up_t&=&u_\dt v_1p_t=u_\dt({\overline{\lambda}}z+1) p_t\\
&=& u_\dt p_t=p_t
\eneq
for all $ t\in [\dt, \|x\|].$

\end{proof}

\begin{thm}\label{Thomp}
Let $A$ be a simple \CA\, with stable rank one. Suppose that $p$ and
$q$ are two $\sigma$-compact open projections of $A$ such that $[p]=[q].$ Then, there
is a precompact support $\{p_{K,n}\}$ of $p,$ and there is a
continuous path of unitaries $\{w(t): t\in [0,1)\}\subset {\tilde
A}$ satisfying the following: $w(0)=1,$ for any $n,$ there is
$t_n\in (0,1)$ such that
\beq\label{Thomp-1}
w(t)p_{K,n}w(t)^*=w(t_n)p_{K,n}w(t_n)^*\tforal t\in [t_n, 1)
\eneq
and $\{w(t_n)^*p_{K,n}w(t_n)\}$ is a precompact support of $q.$
Moreover,
\beq\label{Thomp-1+}
w(t)p_{K,n}=\lambda(t)w(t_n)p_{K,n}
\eneq
for some $\lambda(t)\in \C$ if $t\in [t_n, 1).$
\end{thm}

\begin{proof}
Suppose that $[p]=[q].$ If $p$ is a projection in $A,$ so is $q.$
Then the result follows from a theorem of L. G. Brown (Theorem 1 of
\cite{Br2}).

So, one now assumes that neither $p$ nor $q$ are projections in $A.$
Let $a\in Her(p)$ be a strictly positive element. Let $p_{K,n}$  be
the spectral projection of $a$ associated with $(1/2^{n+1}, \|a\|].$
Then $\{p_{K,n}\}$ is a precompact support for $p.$
 Suppose that $w\in
A^{**}$ such that
\beq\label{Thomp-3}
w^*w=p,\,ww^*=q\andeqn wbw^*\in Her(q)\tforal b\in Her(p).
\eneq
Put $x=wa^{1/2}.$ Then $xx^*=waw^*$ is a strictly positive element
of $Her(q).$

 Put $s_1=1/\sqrt{2},$ $s_n=1/2^{n-1},$ $n=1,2,....$ Since one
 assumes that $p$ is not a projection in $A,$ $0$ is a limit point of
 ${\rm sp}(a^{1/2}).$
Let $p_{s_n}$ be the open spectral projection of $|x|=a^{1/2}$
associated with $(s_n, \|x\|].$ Then $p_{s_n}=p_{K, n},$
$n=1,2,....$
Let $t_n=s_n-s_n/16^n$ and let $p_{t_n}$ be the open spectral
projection of $|x|$ associated with $(t_n, \|x\|],$ $n=1,2,....$

It follows from \ref{PK1} (see also \ref{Ped}) that there is a
unitary $u_n\in {\tilde A}$ with $[u_n]=0$ in $K_1(A)$ such that
\beq\label{Thomp-5}
u_mp_{t_n}u_m^*=q_{t_n}\andeqn u_mp_{t}=u_np_t, {\rm if}\,\,\,m\ge n, t\ge t_n, n, m=1,2,...,
\eneq
where $q_{t_n}$ is the open spectral projection of $|x^*|$
associated with $(t_n, \|x^*\|].$

Denote by $q_{s_n}$ the spectral projection of $|x^*|$ associated with $(s_n, \|x^*\|],$ $n=1,2,....$
Since  $[u_1]=0$ in $K_1(A)$ and $A$ has stable rank one, by a
result of Rieffel (\cite{Rf}),  $u_1\in U_0({\tilde A}).$ Therefore
there is a continuous path of unitaries $\{w(t): t\in [0, t_1]\}\subset {\tilde A}$
($0<t_1<1$) such that
\beq\label{Thomp-9}
w(0)=1, \andeqn w(t_1)=u_1.
\eneq
On also has that
\beq\label{Thomp-12}
u_2u_1^*q_{t_1}&=&u_2u_1^*(u_1p_{t_1}u_1^*)\\\label{Thomp-12+}
 &=&
u_2p_{t_1}u_1^*=u_1p_{t_1}u_1^*=q_{t_1}\\\label{Thomp-12++}
&=& u_1p_{t_1}u_1^*=u_2p_{t_1}u_1^*=u_2p_{t_1}u_2^*(u_2u_1^*)=q_{t_1}(u_2u_1^*).
\eneq
Moreover,
\beq\label{Thomp-13}
u_2u_1^*q_{t_2}(u_2u_1^*)^*=u_2u_1^*(u_1p_{t_2}u_1^*)u_1u_2^*= u_2p_{t_2}u_2^*=q_{t_2}.
\eneq
Let $e_1=p_{(0,(t_1+s_1)/2)}$ the spectral projection of $|x|$ corresponding to $(0, (t_1+s_1)/2)$ and let $C=Her(e_1).$ By
(\ref{Thomp-12+}) and (\ref{Thomp-12++}), one may view $u_2u_1^*$ as a unitary in ${\widetilde
C}.$ Since $C$ has stable rank one and $[u_2u_1^*]=0$ in $K_1(A),$
one obtains a continuous path of unitaries $\{W(t): [t_1,
t_2]\}\subset {\widetilde C} $ ($t_1<t_2<1$) such that
\beq\label{Thomp-15}
W(t_1)=1\andeqn W(t_2)=u_2u_1^*.
\eneq
Note that
\beq\label{Thomp-15+}
q_{s_1}e_1=e_1q_{s_1}=0.
\eneq
 $W(t)$ may be viewed as unitaries in ${\tilde A}.$
Moreover, by (\ref{Thomp-15+}),
\beq\label{Thomp-16-}
W(t)q_{s_1}=\lambda(t)q_{s_1}=q_{s_1}W(t)
\eneq
for some  $\lambda(t)\in \C$ for all $t\in [t_1, t_2].$ Now extend
$w(t)$ from a continuous path from $[0, t_1]$ to a continuous path
from $[0, t_2]$ by defining
\beq\label{Thomp-16}
w(t)=W(t)w(t_1)\tforal t\in [t_1, t_2].
\eneq
Note that
\beq\label{Thomp-17}
w(t_2)p_{s_2}w(t_2)^*=q_{s_2}\andeqn w(t)p_{s_1}w(t)^*=q_{s_1}
\eneq
for all $t\in [t_1, t_2].$ Moreover, by (\ref{Thomp-16-}),
\beq\label{Thomp-17+}
w(t)p_{s_1}&=&W(t)w(t_1)p_{s_1}=W(t)w(t_1)p_{s_1}w(t_1)^*w(t_1)\\
&=& W(t)q_{s_1}w(t_1)=\lambda(t)q_{s_1}w(t_1)\\
&=&\lambda(t)w(t_1)p_{s_1}w(t)^*w(t_1)=\lambda(t)w(t_1)p_{s_1}.
\eneq
Furthermore,
\beq\label{Thomp-18}
w(t_2)=u_2u_1^*w(t_1)=u_2u_1^*u_1=u_2.
\eneq
One also has that
\beq\label{Thomp-19}
u_3w(t_2)^*q_{t_2}&=&u_3u_2^*q_{t_2}u_2u_2^*=u_3p_{t_2}u_2^*=u_2p_{t_2}u_2^*=q_{t_2}\\\label{Thomp-19+}
&=&u_2p_{t_2}u_2^*=u_3p_{t_3}u_3^*(u_3u_2^*)=q_{t_2}(u_3w(t_2)^*),\andeqn
\eneq
\beq\label{Thomp-20}
u_3w(t_2)^*q_{t_3}w(t)u_3^*&=&u_3w(t_2)^*(w(t_2)p_{t_3}w(t_2)^*)w(t_2)u_3^*\\
&=& u_3p_{t_3}u_3^*=q_{t_3}.
\eneq
Therefore, by induction, one obtains a continuous path of unitaries
$\{w(t): t\in [0,1)\}$ of ${\tilde A}$ such that
\beq\label{Thomp-21}
w(0)=1, w(t_n)p_{s_n}w(t_n)^*=q_{s_n}\andeqn\\
w(t)p_{s_n}w(t)^*=q_{s_n} \tforal t\in [t_n, 1).
\eneq
Moreover,
\beq\label{Thomp-22}
w(t)p_{s_n}=\lambda(t)w(t_n)p_{s_n}\tforal t\in [t_n, 1)
\eneq
for some $\lambda(t)\in \C.$

\end{proof}

\begin{cor}\label{compact}
Let $A$ be a simple \CA\, with stable rank one and let $a, \, b\in A_+.$
Suppose that $[p_a]\le [p_b].$ Then, for any $c\in Her(a)_+$ which is compactly contained in $p_a,$
there exists a continuous path of unitaries $\{w(t): t\in [0,1]\}$ such that
$w(0)=1$ and $w(1)^*p_cw(1)$ is compactly contained in $p_b.$

\end{cor}
\begin{proof}

Suppose that $c\in Her(p_a)_+$ which is compactly contained in $Her(p_a)$  in the sense that there is $d\in Her(p_a)_+$ such that
$cd=c.$ Then, for any $\ep>0,$ there is an integer $n\ge 1$ such that
\beq\label{compact-2}
\|f_{1/n}(a)d-d\|<\ep/2.
\eneq
It follows from \cite{R1}  that
\beq\label{compact-3}
f_{\ep}(d)\lessapprox  f_{1/n}(a).
\eneq
Since $A$ has stable rank one, then there exists $v\in {\tilde A}$ such that
$v^*f_{\ep}(d)v\le f_{1/n}(a).$
By \ref{cut},
\beq\label{compact-4}
cf_{\ep}(d)=c.
\eneq
Let $p_c$ be the range projection of $c$ in $A^{**}.$ Then, $v^*p_cv\le p_{K, n}.$
Since   $A$  is simple and has stable rank one, there is $v_0\in {\widetilde{Her(p_c)}}$ such that
\beq\label{compact_5}
[v_0]=[v^*]\,\,\, {\rm in}\,\,\, K_1(A).
\eneq
One may also view $v_0$ as a unitary in ${\tilde A}.$  There is a continuous path of unitaries
$\{w_0(t): t\in [0,1]\}\subset {\tilde A}$ such that
\beq\label{compact-6}
w_0(0)=1, w_0(1)=v_0v.
\eneq

Then
\beq\label{compact-7}
w_0(1)^*p_cw_0(1)=v^*p_cv\le p_{K,n}.
\eneq
Now the  lemma follows  from \ref{Thomp}.

\end{proof}

\begin{prop}\label{Phomp}
The converse of Theorem \ref{Thomp} also holds in the following
sense. Let $A$ be a \CA\, and let $p$ and $q$ be two
$\sigma$-compact open projections of $A.$ Suppose that there is a
continuous path of unitaries $\{w(t): t\in [0,1)\}\subset M(A)$ such
that, $w(0)=1,$ for any $n\ge 1,$ there is $t_n\in (0,1)$ such that
\beq\label{Phomp-1}
w(t)p_{K,n}w(t)^*=w(t_n)p_{K,n}w(t_n)^*\tforal t\in [t_n, 1),
\eneq
where $\{p_{K,n}\}$ is a precompact support for $p$ and
$\{w(t_n)^*p_{K,n}w(t_n)\}$ is a precompact support for $q.$
Moreover,
\beq\label{Phomp-1+}
w(t)p_{K,n}=\lambda(t)w(t_n)p_{K,n}\tforal t\in [t_n, 1)
\eneq
for some $\lambda(t)\in \C.$

Then $[p]=[q].$

\end{prop}

\begin{proof}
One may assume that $t_{n+1}>t_n,$ $n=1,2,....$ Suppose that $a\in
Her(p)$ is a strictly positive element and suppose $s_n\in (0,
\|a\|]$ such that $s_n\searrow 0$ such that $p_{K,n}$ is the
spectral projection of $a$ corresponding to $(s_n, \|a\|],$
$n=1,2,....$ One defines, with $p_{K,0}=0,$
\beq\label{Phomp-2}
v=\sum_{n=1}^{\infty} w(t_n)(p_{K,n}-p_{K,n-1}).
\eneq
One checks that, for $b\in Her(p),$
\beq\label{Phomp-3}
\lim_{n\to\infty}\|\sum_{k=n}^{n+m} w(t_n)(p_{K,n}-p_{K,n-1})b\|=0.
\eneq
It follows that $v\in A^{**}.$
One also checks that
\beq\label{Phomp-2+}
&&\hspace{-1in}\sum_{n=1}^{m}w(t_n)(p_{K,n}-p_{K,n-1})w(t_m)^*\\
&=&\sum_{n=1}^m\overline{\lambda(t_m)}w(t_m)(p_{K,m}-p_{K,n-1})w(t_m)^*\\
&=&\overline{\lambda(t_m)}w(t_m)p_{K,m}w(t_m)^*=\overline{\lambda(t_m)}q_{K,m}.
\eneq

Let $a_n=f_{t_n}(a),$ $n=1,2,...$ Then, if $m>n+1,$
\beq
va_nv^*&=&
{\overline{\lambda(t_m)}}q_{K,m}w(t_m)a_nw(t_m)^*q_{K,m}{\overline{\lambda(t_m)}}\\\label{Phomp-4}
&=& w(t_m)a_nw(t_m)^*\in Her(q).
\eneq
Let $\ep>0$ and $b\in Her(p).$
There is $n\ge 1$ such that
\beq\label{Phomp-5}
\|b-a_mba_m\|<\ep\tforal m\ge n.
\eneq
Then,
\beq\label{Phomp-6}
\|vbv^*-va_mba_mv^*\|<\ep.
\eneq
But, by (\ref{Phomp-4}),  $va_mba_mv^*\in Her(p).$ This implies that $vbv^*\in Her(p).$
Furthermore,
\beq\label{Phomp-7}
v^*pv=q.
\eneq

\end{proof}

\begin{rem}\label{Rhomp}
If $A$ is a unital \CA\, and $p, q\in A$ are two projections which are homotopy, i.e.,
there is a projection $P\in C([0,1], A)$ such that $P(0)=p$ and $P(1)=q.$ 
Then (see, for example, Lemma 2.6.6 of \cite{Lnbk}), there is a unitary 
$U\in C([0,1], A)$ such that $U(0)=1$ and $U(t)^*pU(t)=p(t)$ for all $t\in [0,1].$ 
In \ref{Thomp} and in \ref{Phomp}, for each $n\ge 1,$  there is a
continuous path $$\{p(t)=w(t)^* p_{K,n}w(t): t\in [0, 1)\}$$ of open
projections such that
$$
p(0)=p_{K,n}\andeqn p(t)=q_{K,n}\tforal t\in [t_n,1).
$$
In particular, there is a unitary, namely $u(t_n)\in {\tilde A}$ or
in $M(A)$ in \ref{Phomp}), such that
\beq\label{Rhomp-1}
u(t_n)^*p_{K,n}u(t_n)=q_{K,n},\,\,\,n=1,2,....
\eneq

In general, however, if $[p]=[q]$ in the sense of Cuntz, there may
not be any unitary path $\{w(t): t\in [0,1]\}$ for which
$w(0)^*pw(0)=p$ and $w(1)^*pw(1)=q,$ as one can see from the
following.

\end{rem}

\begin{prop}\label{Countex}
Let $A$ be a non-unital and $\sigma$-unital non-elementary simple
\CA\, with (SP).  Then there are two $\sigma$-compact open
projections $p$ and $q$ of $A$ such that $[p]=[q]$ but there are no
unitary $u\in M(A)$ such that $u^*pu=q.$

\end{prop}

\begin{proof}
Let $a\in A_+$ be a strictly positive element. Then  $0$ is a limit
point of ${\rm sp}(a).$
Thus $A$ admits an approximate identity $\{e_n\}$ such that
$e_{n+1}e_n=e_n,$ $n=1,2,....$ One may further assume, without loss
of generality, by passing to a subsequence if necessary,
 there are nonzero positive element $b_n\in
 Her(e_{2(n+1)}-e_{2n})$ with $\|b_n\|=1,$ $n=1,2,....$ In particular,
 \beq\label{Countex-0}
 b_ib_j=0\,\,\,{\rm if}\,\,\, i\not=j.
 \eneq
 On the other hand, since $\overline{e_1Ae_1}$ is a non-elementary
simple \CA, by a result of Akemann and Shultz \cite{AS}, there are
mutually orthogonal non-zero positive elements $c_1,
c_2,...,c_n,...$ in $\overline{e_1Ae_1}.$ Since $A$ has (SP), there
are non-zero projections $d_n'\in Her(c_n),$ $n=1,2,....$

By a result of Cuntz (see (2) of Lemma 3.5.6 of \cite{Lnbk}, for
example), there are partial isometries  $x_1, x_2,...,x_n,...\in A$
such that
\beq\label{Countex-1}
x_i^*x_i\in Her(c_i)\andeqn x_ix_i^*\in Her(b_i),
\eneq
where $x_i^*x_i$ and $x_ix_i^*$ are non-zero projections,
$i=1,2,....$

Put  $d_n=x_nx_n^*$ and $f_n=x_n^*x_n,$ $n=1,2,....$ Define
\beq\label{Countex-2}
b=\sum_{n=1}^{\infty}{y_n\over{n^2}}\andeqn c=\sum_{n=1}^{\infty}
{z_n\over{n^2}}.
\eneq
Then $b,c\in A.$  Define $x=\sum_{n=1}^n {x_n\over{n}}.$ Then
\beq\label{Countex-3}
x^*x=b\andeqn xx^*=c.
\eneq
Let $p=p_c,$ the range projection of $c$ in $A^{**}$ and let
$q=p_b,$ the range projection of $b$ in $A^{**}.$ Then, by
(\ref{Countex-3}), $[p]=[q].$ Moreover $ce_2=c.$ So $c$ is compact.
Furthermore, since $f_n\le (e_{2(n+1)}-e_{2n}),$ $n=1,2,...,$
\beq\label{Countex-3+}
q=\sum_{n=1}^{\infty}f_n,
\eneq
where the sum converges in the strict topology.
It follows
that $q\in M(A).$

Now suppose that there were a unitary $u\in M(A)$ such that
$u^*qu=p.$ Therefore  $u^*qu\in M(A).$ However   $p\not\in M(A).$
Otherwise $pe_2=p$ implies that $p\in A.$ But $p\not\in A.$ So there is no unitary $u\in M(A)$ for which
$u^*qu=p.$

\end{proof}

\begin{rem}
In the proof of \ref{Countex}, one notes that $p$ is precompact and
is compactly contained in $p_{e_2}.$ However, $q$ is not precompact
and is not compactly contained in any $\sigma$-compact open
projection of $A.$ In fact, if $q\le a$ for some $a\in A_+, $ then
$q\in A$ since $q\in M(A).$ One concludes that precompactness is not
invariant under the Cuntz relation.

\end{rem}

Finally, to end this section,  one has the following:

\begin{prop}\label{Fhomp}
Let $A$ and $B$ be two separable \CA s, and let $\phi_0, \phi_1:
A\to B$ be  two \hm s. Suppose that there is a \hm\, $H:  A\to
C([0,1], B)$ such that $\pi_0\circ H=\phi_0$ and $\pi_1\circ
H=\phi_1,$ where $\pi_t: C([0,1], B)\to B$ is the point-evaluation
at the point $t\in [0,1].$   Suppose also that $H$ extends to a
(sequentially) normal \hm\, $H': A^{**}\to C([0,1], B^{**})$ in the
sense that if $\{a_n\}\subset A_{s.a}$ is a increasing bounded
sequence with upper bound $x\in A^{**},$ then $\{H(a_n)\}$ has the
upper bound $H'(x).$ Then $\phi_0$ and $\phi_1$ induce the same
\hm\, on the Cuntz semigroups $W(A)$ and $Co(A).$
\end{prop}

\begin{proof}
Let $p$ be an open projection of $A.$ Since $A$ is separable, $p$ is
$\sigma$-compact. Let $a\in A_+$ such that $p$ is the range
projection of $A.$ Then $\{a^{1/n}\}$ is increasing and has the upper
bound $p.$

Put $p(t)=\pi_t\circ H'(p),$ $t\in [0,1].$ It is a norm continuous
path of $\sigma$-compact open projections.  It follows from
\ref{Hp1} that $[p(0)]=[p(1)].$ In other words,
\beq\label{Fhomp-1}
[\phi_0(p)]=[\phi_1(p)].
\eneq

\end{proof}

\section{Hilbert Modules}

From the definition (see \ref{Dcuntz}), two $\sigma$-compact open projections 
$p$ and $q$ of a \CA\, $A$ are Cuntz equivalent if and only if the corresponding Hilbert 
$A$-modules are isomorphic as Hilbert $A$-modules. When $A$ has stable rank one, 
by \cite{CEI} (see also \ref{TCEI} ), two positive elements $a$ and $b$ in $A$ are Cuntz equivalent 
if and only if the associated Hilbert $A$-modules are isomorphic as Hilbert $A$-modules. 
A question  was mentioned in \cite{CEI} (see line 27 of page 187 of \cite{CEI}) whether $a$ and $b$ are Cuntz equivalent if the Banach $A$-modules $\overline{aA}$ and 
$\overline{bA}$ are isomorphic as Banach modules.  This question will be answered by 
a result of L. G. Brown below (\ref{Brown} and \ref{Rbr}). Recently, related to the Cuntz semigroups,
projective Hilbert modules 
also attract some attention (see \cite{BC}).  In this section, 
these two issues will be discussed.  We begin with the following definition. 

\vspace{0.1in}

\begin{def}\label{Hild}
Let $A$ be a \CA. For an integer $n\ge 1,$ denote by $A^{(n)}$ the
Hilbert $A$-module of orthogonal direct sum of $n$ copies of $A.$ If
$x=(a_1, a_2,...,a_n), y=(b_1,b_2,...,b_n),$ then
$$
<x,y>=\sum_{i=1}^n a_n^*b_n.
$$
Denote by $H_A$ the standard countably generated Hilbert (right)
$A$-module
$$
H_A=\{\{a_n\}: \sum_{n=1}^k a_n^*a_n \,\,\,{\rm converges\,\,\,
in\,\,\, norm}\},
$$
where the inner product is defined by
$$
<\{a_n\}, \{b_n\}>=\sum_{n=1}^{\infty} a_n^*b_n.
$$

Let $H$ be a Hilbert $A$-module. denote by $H^{\sharp}$ the set of
all  bounded $A$-module maps from $H$ to $A.$ If $H_1, H_2$ are
Hilbert $A$-modules, denote by $B(H_1, H_2)$ the space of all
bounded module maps from $H_1$ and $H_2.$   If $T\in B(H_1, H_2),$
denote by $T^*: H_2\to H_1^{\sharp}$  the bounded module maps
defined by
$$
T^*(y)(x)=<Tx, y>\tforal x\in H_1\andeqn y\in H_2.
$$
If $T^*\in B(H_2, H_1),$ one says that $T$ has an adjoint $T^*.$  Denote by
$L(H_1, H_2)$ the set of all bounded $A$-module maps in $B(H_1,
H_2)$ with  adjoints.  Let $H$ be a Hilbert $A$-module. In what
follows, denote $B(H)=B(H,H)$ and $L(H)=L(H, H).$  $B(H)$ is a
Banach algebra and $L(H)$ is a \CA.

Denote by $F(H)$ the linear span of those module maps with the form
$\xi<\zeta, ->,$ where $\xi, \zeta\in H.$ Denote by $K(H)$ the
closure of $F(H).$  $K(H)$ is a \CA.  It follows from a result of
Kasparov (\cite{K}) that $L(H)=M(K(H)),$  the multiplier algebra of
$K(H),$ and, by \cite{Lnbd}, $B(H)=LM(K(H)),$ the left multiplier
algebra of $K(H).$

Two Hilbert $A$-modules are said to be unitarily equivalent, or isomorphic,  if there
is an invertible map $U\in B(H_1, H_2)$ such that
$$
<U(x_1), U(x_2)>=<x_1, x_2>\tforal x_1, x_2\in H_1.
$$
\end{def}

The  following result of L. G. Brown
becomes quite useful and answers the question in \cite{CEI} mentioned above.

\begin{thm}\label{Brown}  {\rm (Theorem 3.2 of
\cite{Br1} and Theorem 2.2 of \cite{Lninj})} Let $H_1$ and $H_2$ be
two countably generated Hilbert modules over a \CA\, $A.$ Suppose
that there is $T\in B(H_1, H_2)$ which is one-to-one and has dense
range. Then $H_1$ and $H_2$ are unitarily equivalent.
\end{thm}

\begin{rem}\label{Rbr}
However, it is also worth to note that the above statement fails
when $H_1$ and $H_2$ are not countably generated. See Example 2.3 of
\cite{Lninj}.
\end{rem}

 As a consequence, one has the following.

\begin{prop}\label{pq}
Let $A$ be a \CA\, and let  $a, b\in A_+.$   Suppose that
$H_1=\overline{aA}$ and $H_2=\overline{bA}.$ Then $[p_a]\le [p_b]$
(or equivalently, $a\lessapprox b$) if and only if there is $T\in
B(H_2, H_1)$ whose range is  dense in $H_1.$

\end{prop}

\begin{proof}
Suppose that $[p_a]\le [p_b],$ i.e., there is a partial isometry $v\in A^{**}$ such that
\beq\label{pq-1-1}
v^*p_av\le p_b\andeqn v^*xv\in Her(b)\tforal x\in Her(a).
\eneq
Thus $v^*H_1\subset H_2.$ Put $H_3=\overline{v^*H_1}$ and  $c=v^*av.$ Then $c\in  K(H_3).$
It follows from Lemma 2.13 of \cite{Lninj} that one may view $K(H_3)$ as a hereditary \SCA\, of $K(H_1).$
Thus $T=vc$ defines a bounded module map in $B(H_2, H_1).$  Note $T=av$ and $\overline{vH_2}=H_1.$
It follows that $T$ has the dense range.

Now one assumes that there is $T\in B(H_2, H_1)$ whose range is dense in $H_1.$
One may identify $T$ with an element in $LM(Her(b), Her(a)).$ Let
$x=(Tb)^*Tb.$ Then $x\in Her(b).$ Let $H_4=\overline{xA}.$ Then $T$
is one-to-one on $H_4$ and has dense range. It follows from
\ref{Brown} that $H_4$ and $H_1$ are unitarily equivalent which
provides a partial isometry $v\in A^{**}$ such that
\beq\label{pq-1}
v\overline{aA}=\overline{xA}\andeqn\hspace{-0.05in},\,\,\text{for}\,\,\,
\xi\in \overline{aA},\,\,\, v\xi=0\,\,\,\text{if and only if}\,\,\,
\xi=0
\eneq
Let $r=p_x.$ Then
\beq\label{pq-2}
vp_av^*=r\le p_b\andeqn v\xi v^*\in Her(x)\subset Her(b) \tforal \xi\in Her(a).
\eneq

\end{proof}

%
%
%

Now we turn to the projectivity of Hilbert modules.

\begin{thm}\label{1proj}
Let $A$ be a \CA. Then every Hilbert $A$ module is projective (with
bounded module maps with adjoints as morphisms) in the following
sense: Let $H$ be a Hilbert $A$-module.

{\rm (1)}  Suppose that $H_1$ is another Hilbert $A$-module and
suppose that $\phi\in L(H_1, H)$ is a surjective.  Then there is
$\psi\in L(H, H_1)$ such that
\beq\label{1proj-1}
\phi\circ \psi={\rm id}_{H};
\eneq

{\rm (2)} Suppose that $H_2$ and $H_3$ are Hilbert modules and
suppose that $\phi_1\in L(H_2, H_3)$ is surjective. Suppose also
that $\phi_2\in L(H, H_3).$ Then there exists  $\psi\in L(H, H_2)$
such that
\beq\label{1proj-2}
\phi_1\circ \psi=\phi_2.
\eneq

\end{thm}

\begin{proof}
For (1), one first notes that $\phi$ has closed range.  Define $T:
H_1\oplus H\to H_1\oplus H$ by $T(h_1\oplus h)=0\oplus T(h_1)$ for
$h_1\in H_1$ and $h\in H.$ Then $T\in L(H_1\oplus H)=M(K(H_1\oplus
H)).$  It follows from Lemma 2.4 of \cite{Lninj} that
\beq\label{1proj-3}
H_1\oplus H={\rm ker}T\oplus |T|(H_1\oplus H).
\eneq
Let $T=V|T|$ be the polar decomposition in $(K(H_1\oplus H))^{**}.$
Note that the proof of Lemma 2.4 of \cite{Lninj} shows that $0$ is an
isolated point of $|T|$ or $|T|$ is invertible. So the same holds for
$(TT^*)^{1/2}.$  Let $S=(TT^*)^{-1},$ where the inverse is taken in
the hereditary \SCA\, $L(H)\subset L(H_1\oplus H).$ Since $T$ is
surjective,
\beq\label{1proj-4}
|TT^*|H=H.
\eneq
Moreover,
\beq\label{1proj-5}
L_1=V^*(TT^*)^{-1/2}=V^*(TT^*)^{1/2}S\in L(H_1\oplus H).
\eneq

One then checks that
\beq\label{1proj-7}
TL_1= V|T|V^*(TT^*)^{-1/2}=P,
\eneq
where $P$ is the range projection of $(TT^*)^{1/2}$ which gives the
identity of $H.$ One then defines $\psi$ by $L_1.$ Thus $\phi\circ
\psi={\rm id}_H.$

For (2), one applies (1). Since  $\phi_1$ is surjective, by (2),
there is $\phi_3\in L(H_3, H_2)$ such that
\beq\label{1proj-8}
\phi_1\phi_3={\rm id}_{H_3}.
\eneq
Define $\psi=\phi_1\circ \psi_3\circ \phi_2.$

\end{proof}

\begin{rem}\label{R1}
A discussion  about injective Hilbert modules can be found in
\cite{Lninj}.  It was shown that, for example, a Hilbert $A$-module
$H$ is injective (with bounded module maps with adjoints as
morphisms) if and only if it is  orthogonally  complementary
(Theorem 2.14 of \cite{Lninj}). For a full countably generated
Hilbert module, it is injective (with bounded morphisms with
adjoints as morphisms) if and only if $L(H)=B(H)$ (see 2.9 and 2.19
of \cite{Lninj}).

Let $A$ be a \CA. One may consider the category of Hilbert
$A$-modules with bounded $A$-module maps as morphisms. A discussion
on the question which Hilbert $A$-modules are injective in this
category was given in \cite{Lninj}. It seems that  question which
Hilbert $A$-modules are projective in this category is much more
difficult. Consider a Hilbert $A$-module $H=\xi A$ which is singly
algebraically generated. Let $H_1$ be another Hilbert $A$-module and
$T\in B(H_1, H)$ is surjective. Suppose that $x\in H_1$ such that
$T(x)=\xi.$ It would be most natural to define $S: H\to H_1$ by
$S(\xi)=x$ which gives $TS(y)=y$ for all $y\in H.$ The trouble is
that it is not clear why $S$ should be bounded.

Noticing  the  difference between algebraically projective
$A$-modules and projective Hilbert $A$-modules (with bounded module
maps as morphisms),
 the following two propositions may not seem entirely trivial.
The first one is certainly known. 
After this note was first posted, Leonel Robert informed the author 
that, using Proposition \ref{1proj} above, he has a proof that the converse of the following also 
holds, i.e., if $H$ is algebraically finitely generated, then $K(H)$ has an identity. 
\end{rem}

\begin{prop}\label{PH}
Let $A$ be a \CA\, and let $H$ be a Hilbert $A$-module. Suppose that
$1_{H}\in K(H).$ Then $H$ is algebraically finitely generated.
\end{prop}

\begin{proof}
Let $F(H)$ be the linear span of  rank one module maps of the form
$\xi<\zeta, ->$ for $\xi, \zeta\in  H.$ Then $F(H)$ is dense in $K(H).$ There is $T\in
F(H)$ such that
\beq\label{PH-1}
\|1_H-T\|<1/4,
\eneq
One may assume that $\|T\|\le 1.$ Thus
\beq\label{PH-2}
\|1_H-T^*T\|<1/2.
\eneq
It follows that   $0\le T^*T\le 1_H$ and $T^*T$ is invertible. Note that
$T^*T\in F(H).$ Therefore there are $\xi_1, \xi_2,...,\xi_n,
\zeta_1,\zeta_2,...,\zeta_n\in H$ such that
\beq\label{PH-3}
T^*T(\xi) =\sum_{j=1}^n \xi_j<\zeta_j, \xi>\tforal \xi\in H.
\eneq
But $T^*TH=H.$ This implies that $\sum_{j=1}^n\xi_jA=H.$

\end{proof}


\begin{prop}\label{FHP}
Let $A$ be a \CA\, and let $H$ be a Hilbert $A$-module for which
$K(H)$ has an identity. Then $H$ is projective Hilbert $A$-module
(with bounded module maps as morphisms).
\end{prop}

\begin{proof}
One first assumes that $A$ has an identity. From \ref{PH}, $H$ is
finitely generated. Therefore, a theorem of Kasparov shows that
$H=PH_A$ for some projection $P\in L(H_A).$ The fact that $1_H\in
K(H)$ implies that $P\in K(H_A).$ Therefore there is an integer
$N\ge 1$ and a projection $P_1\in M_N(A)$ such that $PH$ is
unitarily equivalent to $P_1H_A.$  In other words, one may assume
that $H$ is a direct summand of $A^{(N)}.$ Suppose that $H_1$ and
$H_2$ are  two Hilbert $A$-modules and suppose that $S\in B(H_1,
H_2)$ is surjective and suppose that  $\phi: H\to H_2$ is a bounded
module map. Since $H$ is a direct summand of $A^{(N)},$ there is a
partial isometry $V\in L(H, A^{(N)})$ such that $P_1V={\rm id}_H.$
Let $T=\phi\circ P_1.$ Denote by $e_i$ the vector in the $i$th copy
of $A$ given by $1_A.$ Choose $g_1, g_2,...,g_n\in H_1$ such that
$Sg_i=Te_i,$ $i=1,2,...,n.$ Define $L: K(A^{(N)}\oplus H_1)$ by
\beq\label{FHP-1}
L(h\oplus h_1)=\sum_{j=1}^N g_i<e_i, h>\tforal h\in H\andeqn h_1\in
H_1.
\eneq
Define $L_1=L|_H.$ For  $h=\sum_{j=1}^Ne_ia_i,$ where $a_i\in A, $
one has
\beq\label{FHP-2}
SL_1(h)&=& S(\sum_{j=1}^N g_i<e_i, h>)
=\sum_{j=1}^N Sg_i<e_i,h>\\
&=&\sum_{j=1}^N Te_i<e_i, e_i>a_i=\sum_{j=1}^N Te_ia_i\\
&=&T(h).
\eneq
Define $L_2\in B(H, H_1)$ by $L_2=L_1\circ V.$ Then
\beq\label{FHP-3}
SL_2=SL_1\circ V=T\circ V=\phi\circ P_1\circ V=\phi.
\eneq

Moreover, if $S_1\in B(H_1, H)$ is a surjective map, consider the
following diagram:
$$
\begin{array}{ccccc}
&& H &\\
& & \hspace{0.1in}\downarrow_{{\rm id}_H} &\\
H_1 & \rightarrow_{S_1} & H  &\to 0
\end{array}
$$
From (ii), there is a bounded module map $L: H\to H_1$ such that
\beq\label{FHP-4}
S_1L={\rm id}_{H}.
\eneq

For general case, one may consider $H$ as a Hilbert  ${\tilde
A}$-module.

\end{proof}

\begin{rem}
The fact that $<e_i, e_i>=1_A$ is crucial in the proof. It should be
noted that, when $A$ is not unital,  the above argument does not
imply  that $A^{(n)}$ is projective (with bounded module maps as
morphisms).

\end{rem}

\begin{cor}\label{Ansub}
Let $A$ be a unital \CA\, and let $H$ be a Hilbert $A$-module.
Suppose that there is an integer $n\ge 1$ and a surjective map $S\in
B(A^{(n)}, H).$ Then $H$ is projective (with bounded module maps as
morphisms)
\end{cor}

\begin{proof}
Let $H_1$ and $H_2$ be two Hilbert $A$-modules and let $\phi\in
B(H_1, H_2)$ which is surjective. Suppose that $\psi\in B(H, H_2).$

Since $A^{(n)}$ is self-dual, $S^*$  must map $H$ into $A^{(n)}.$
In other words,  $S\in L(A^{(n)}, H).$ By  \ref{1proj}, there exists
$T\in L(H, A^{(n)})$ such that
\beq\label{Ansub-1}
ST={\rm id}_H.
\eneq
Let $\phi_1\in B(A^{(n)}, H_2)$ be defined by
\beq\label{Ansub-2}
\phi_1=\phi\circ S.
\eneq
Then, by \ref{FHP}, $A^{(n)}$ is projective. There is  $L\in
B(A^{(n)}, H_1)$ such that
\beq\label{Ansub-3}
\phi\circ L=\phi_1.
\eneq
Define $\phi_2=L\circ T.$ Then $\phi_2\in B(H, H_1).$ Moreover,
\beq\label{Ansub-4}
\phi\circ \phi_2 &=& \phi\circ L\circ T=\phi\circ S\circ T=\phi.
\eneq
Hence $H$ is projective (with bounded module maps as morphisms).

\end{proof}

There are projective Hilbert modules (with bounded module maps as
morphisms) for which $K(H)$ is not unital.

\begin{thm}\label{LM=M}
Let $A$ be a separable \CA\, such that $LM(A\otimes {\cal
K})=M(A\otimes {\cal K}).$ Then every countably generated Hilbert
$A$-module is projective (with bounded module maps as morphisms)
\end{thm}

One needs the following lemma which the author could not locate a
reference.

\begin{lem}\label{openmap}
Let $X$ be a Banach space and let $H$ be a separable Banach space.
Suppose that $T: X\to H$ is a surjective bounded linear map.  Then
there is a separable subspace $Y\subset X$ such that $TX=H.$
\end{lem}

\begin{proof}
Note that the Open Mapping Theorem applies here. From the open
mapping theorem (or a proof of it), there is $\dt>0$ for which
$T(B(0, a))$ is dense in $O(0, a\dt)$ for any $a>0,$  where $B(0,
a)=\{x\in X: \|x\|\le a\}$ and $O(0, b)=\{h\in H: \|h\|<b\}.$ For
each rational number $r>0,$ since $H$ is separable, one may find a
countable set $E_r\subset B(0, r)$ such that $T(E_r)$ is dense in
$O(0,r\dt).$ Let $Y$ be the closed subspace generated by $\cup_{r\in
\Q_+} E_r.$

Let $d=\dt/2$ and let $y_0\in O(0, d).$ Then $T(Y\cap B(0, 1/2))$ is
dense in $O(0,d).$ Choose $\xi_1\in Y\cap B(0,1/2)$ such that
\beq\label{Opm-1}
\|y_0-T\xi_1\|<\dt/2^2.
\eneq
In particular,
\beq\label{Opm-2}
y_1=y_0-T\xi_1\in O(0, \dt/2^2).
\eneq
Since $T(Y\cap B(0, 1/2^2))$ is dense in $O(0, \dt/2^2),$ one
obtains $\xi_2\in Y\cap B(0, 1/2^2)$ such that
\beq\label{Opm-3}
\|y_1-T\xi_2\|<\dt/2^3.
\eneq
In other words,
\beq\label{Opm-4}
y_2=y_1-T\xi_2=y_0-(T\xi_1+T\xi_2)\in O(0, \dt/2^3).
\eneq
Continuing this process, one obtains  a sequence of elements
$\{\xi_n\}\subset Y$ for which $\xi_n\in B(0, 1/2^n)$ and
\beq\label{Opm-5}
\|y_0-(T\xi_1+T\xi_2+\cdots +T\xi_n)\|<\dt/2^{n+1},\,\,\,n=1,2,....
\eneq
Define $\xi_0=\sum_{n=1}^{\infty}\xi_n.$ Note that the sum converges
in norm and therefore $\xi_0\in Y.$ By the continuity of $T,$
\beq\label{Opm-6}
T\xi_0=y_0.
\eneq
This implies that $T(Y)\supset O(0, d).$  It follows that $T(Y)=H.$
\end{proof}

{\it Proof of Theorem \ref{LM=M}}

Let $H$ be a countably generated Hilbert $A$-module. Suppose that
$H_1$ and $H_2$ are two Hilbert $A$-modules, suppose that
$\phi\in B(H_1, H_2)$  and $\psi\in B(H, H_2).$ Suppose also that
$\phi$ is surjective.

Let $H_3=\overline{\psi(H)}.$ Then $H_3$ is countably generated.
Since $A$ is separable, $H_3$ is also a separable Banach space. By
\ref{openmap}, there is a separable subspace $Y\subset H_1$ such
that $TY=H_3.$ Let $H_4$ be the Hilbert $A$-module generated by $Y.$
Then $H_4$ is countably generated.

Let $H_0=H_A\oplus H_4\oplus H_3\oplus H.$ Then, by a result of
Kasparov (\cite{K}), $H_0\cong H_A.$ Define
\beq\label{LM=M-1}
\Psi(h_0\oplus h_4\oplus h_3\oplus h)=\psi(h_4)\andeqn
\Phi(h_0\oplus h_4\oplus h_3\oplus h)=\phi(h)
\eneq
for all $h_0\in H_A, h_4\in H_4, h_3\in H_3$ and $h\in H.$ Note that
$\Psi$ is  from $H_0$ onto $H_3.$

By the assumption that $LM(A\otimes K)=M(A\otimes K)$ and  by Theorem 1.5 of \cite{Lnbd} and \cite{K}, $\Psi,
\Phi\in L(H_0).$ It follows that $\phi|_{H_4}\in L(H_4, H_3)$ and
$\psi\in L(H, H_3).$ By \ref{1proj},  there exists
$\phi_1\in L(H, H_4)$ such that
\beq\label{LM-M-2}
\phi\circ \phi_1=\psi.
\eneq

\vspace{0.2in}

\begin{lem}\label{app}
Let $A$ be a \CA\, and  let $H$ be a Hilbert $A$-module. Let
$H_0\subset H$ be  a Hilbert $A$-submodule. Suppose that $\{e_\af\}$
is an approximate identity for $K(H_0)$ and suppose that $\xi\in H.$
Then
\beq\label{app-1}
\|\pi(\xi)\|=\lim_\af\|(1-e_\af)(\xi)\|,
\eneq
where $\pi: H\to H/H_0$ is the quotient map.
\end{lem}

\begin{proof}
 Note that
$$
\|\pi(\xi)\|=\inf\{\|\xi+\zeta\|: \zeta\in H_0\}.
$$
It follows from Lemma 2.13 of \cite{Lninj} that $K(H_0)$ may be regarded as a
hereditary \SCA\, of $K(H).$

Let $\ep>0.$ There exists $\zeta\in H_0$ such that
\beq\label{app-2}
\|\pi(\xi)\|\ge \|\xi+\zeta\|-\ep/2.
\eneq
There exists $\af_0$ such that
\beq\label{app-3}
\|(1-e_\af)(\zeta)\|<\ep/4\tforal \af\ge \af_0.
\eneq
Note that $0\le 1-e_\af\le 1$ for all $\af.$ Therefore
\beq\label{app-4}
\|\pi(\xi)\|&\ge&  \|\xi+\zeta\|-\ep/2 \ge
\|(1-e_\af)(\xi+\zeta)\|-\ep/2\\
&\ge &
\|(1-e_\af)(\xi)\|-\|(1-e_\af)(\zeta)\|-\ep/2\\
&\ge & \|(1-e_\af)(\xi)\|-\ep.
\eneq
Let $\ep\to 0,$
\beq\label{app-5}
\|\pi(\xi)\|\ge \|(1-e_\af)(\xi)\|\tforal \af\ge \af_0.
\eneq
It follows that
\beq\label{app-6}
\|\pi(\xi)\|\ge \lim_\af \|(1-e_\af)(\xi)\|.
\eneq
Since $e_\af(\zeta)\in H_0$ for all $\af,$
\beq\label{app-7}
\|\pi(\xi)\|\le \lim_\af \|(1-e_\af)(\xi)\|.
\eneq
The lemma follows from the combination of (\ref{app-6}) and
(\ref{app-7}).

\end{proof}

\begin{rem}\label{R2}

Suppose that $H_1$ and $H$ are Hilbert $A$-modules and $\phi: H_1\to
H$ is a bounded surjective module map.  Let $H_0={\rm ker} \phi.$ It
is a Hilbert submodule of $H_1.$ Let $\pi: H_1\to H_1/H_0$ be the
quotient map. It is a Banach space. There is a bounded linear map
$\phi': H_1/H_0\to H$ such that $\phi'\circ \pi=\phi.$  Since
$\phi'$ is one-to-one and onto, it has an inverse. In what follows
denote by $\phi^{\sim}: H\to H/H_0$ the inverse  which is also
bounded.

 Let
$p$ be the open projection of $K(H_1)$ corresponding  $K(H_0).$ Then
$H/H_0$ may be identified with $(1-p)H$ which can also be made into
a Banach $A$-module.

\end{rem}

\begin{lem}\label{estimate}Let $A$ be a \CA\, and let $H$ be a Hilbert $A$-module. Suppose
that $\xi_1, \xi_2,...,\xi_n\in H$ and $e_1, e_2,...,e_n\in A_+$ are
in the center of $A$  with $0\le e_i\le 1$ ($i=1,2,...,n$) such that
\beq\label{est-1}
e_ie_j=e_je_i=0\,\,\, \text{if}\,\,\,|i-j|\ge 2\tand
\xi_ie_i=\xi_i,\,\,\,i=1,2,...,n.
\eneq
Then, for any $b\in A,$
\beq\label{est-2}
\|\sum_{i=1}^n \xi_ib\|\le 2 \max_{1\le i\le n}\|\xi_i\|\|b\|.
\eneq
\end{lem}

\begin{proof}
Let $F\in H^{\sharp}$ with $\|F\|\le 1.$ Let $p_i$ be the range
projection of $F(\xi_i)^*F(\xi_i)$ in $A^{**},$ $i=1,2,...,n.$ Note
that $p_ip_j=p_jp_i=0$ if $|i-j|\ge 2.$

Define
$$
C_0=\begin{pmatrix}F(\xi_2) & F(\xi_4) &\cdots & F(\xi_{2k})\\
                    0 & 0 &\cdots & 0\\
                    \vdots & \vdots & & \vdots\\
                    0 & 0 & \cdots & 0\end{pmatrix}
                    \andeqn
                    C_1=\begin{pmatrix}F(\xi_1) & F(\xi_3) &\cdots & F(\xi_{2k-1})\\
                    0 & 0 &\cdots & 0\\
                    \vdots & \vdots & & \vdots\\
                    0 & 0 & \cdots & 0\end{pmatrix}
                    $$
                    $$
                    B_0=\begin{pmatrix} p_2b_2 & 0 &\cdots & 0\\
                    p_4b_4 & 0 &\cdots & 0\\
                    \vdots & \vdots & & \vdots\\
                    p_{2k}b_{2k} & 0 & \cdots & 0\end{pmatrix}
                    \andeqn
                    B_1=\begin{pmatrix} p_1b_1 & 0 &\cdots & 0\\
                    p_3b_3 & 0 &\cdots & 0\\
                    \vdots & \vdots & & \vdots\\
                    p_{2k-1}b_{2k-1} & 0 & \cdots & 0\end{pmatrix}.
                    $$
Here if $n$ is even, then $2k=n,$ if $n$ is odd, then $n=2k-1$ and
$\xi_{2k}=0.$
One estimates that
\beq\label{est-3}
\|F(\sum_{i=1}^n \xi_ib)\| &=&\|\sum_{i=1}^n
F(\xi_i)p_ib\|\\
&\le & \|\sum_{i=odd} F(\xi_i)p_ib\|+\|\sum_{i=even}
F(\xi_i)p_ib\|\\
&=&\|C_1B_1\|+\|C_0B_0\|\le
(\|C_1^*C_1\|\|B_1^*B_1\|)^{1/2}+(\|C_0^*C_0\|\|B_0^*B_0\|)^{1/2}\\
&=&(\|C_1C_1^*\|\|B_1^*B_1\|)^{1/2}+(\|C_0C_0^*\|\|B_0^*B_0\|)^{1/2}\\
&\le &(\|\sum_{i=odd}
F(\xi_i)F(\xi_i)^*\|\|\sum_{i=odd}b^*p_ib\|)^{1/2}\\
&&+
(\|\sum_{i=even} F(\xi_i)F(\xi_i)^*\|\|\sum_{i=even}^nb^*p_ib\|)^{1/2}\\
&=& (\|\sum_{i=odd}
e_iF(\xi_i)F(\xi_i)^*e_i\|\|b^*(\sum_{i=odd}^np_i)b\|)^{1/2}\\
&&+(\|\sum_{i=even}
e_iF(\xi_i)F(\xi_i)^*e_i\|\|b^*(\sum_{i=even}^np_i)b\|)^{1/2}\\
&\le & ((\max_{i=odd} \|F(\xi_i)F(\xi_i)^*\|) \|\sum_{i=odd}
e_i\|\|b^*b\|)^{1/2}\\
&&+ ((\max_{i=even} \|F(\xi_i)F(\xi_i)^*\|) \|\sum_{i=even}
e_i\|\|b^*b\|)^{1/2}\\
&\le & 2\max_{1\le i\le n} \|\xi_i\| (\|b^*b\|)^{1/2}\\
\eneq
It follows that
\beq\label{est-4}
\|\sum_{i=1}^n \xi_ib_i\|\le 2\max_{1\le i\le n} \|\xi_i\| \|b\|
\eneq
\end{proof}

\begin{rem}
In the lemma above, if $e_1,e_2,...,e_n$ are mutually orthogonal,
then the number $2$ in (\ref{est-2}) can be replaced by $1.$

\end{rem}

\begin{df}
Let $A$ be a \CA. An approximate identity $\{e_n\}$ is said to be a {\it sequential central} approximate identity,
if $\{e_n\}$ is a sequence and each $e_n$ is in the center of $A.$ 
\end{df}

\begin{thm}\label{Tproj}
Let $A$ be a unital \CA, let  $a\in A\setminus \{0\}$ and let
$H=\overline{aA}.$ Suppose that $K(H)$ has a sequential central
approximate identity. Then $H$ is a projective Hilbert $A$-module
(with bounded module maps as morphisms).

Moreover, if $H_1$ and $H_2$ are two Hilbert $A$-modules,  $\phi\in
B(H_1, H_2)$ is surjective and if $\psi\in B(H, H_2).$ Then, for any
$\ep>0,$ there exists $T\in B(H, H_1)$ with
$$
\|T\|\le 2\|\phi^{\sim}\circ \psi\|+\ep
$$
such that
$$
\phi\circ T=\psi.
$$

In the case that $K(H)$ admits a central approximate identity
consisting of a sequence of projections, one can choose $T\in B(H,
H_2)$ such that
$$
\|T\|\le \|\phi^{\sim}\circ \psi\|+\ep.
$$
\end{thm}

\begin{proof}
Since $\overline{aAa}$ has a sequential central approximate
identity, $\overline{aAa}$ contains a strictly positive element $x$
which is in the center. One may assume that  $a=x$ and ${\rm
sp}(a)=[0,1].$ Let $f_n\in C_0((0,1])$ be such that $0\le f_n\le 1,$
$f_n(t)=1$ if $t\in [1/2^n, 1],$ $f_n(t)=0$ if $t\in [0,3/2^{n+2}]$
and $f(t)$ is linear in $[3/2^{n+2}, 1/2^n],$ $n=1,2,...,$ and let
$g_n\in C_0((0,1])$ be such that $0\le g_n\le 1,$ $g_n(t)=1$ if
$t\in [1/2^{n+2}-1/2n2^{n+2}, 1/2^n+1/2n2^{n+2}],$ $g_n(t)=0$ if
$t\not\in [1/2^{n+2}-1/n2^{n+2}, 1/2^n+1/n2^{n+2}]$ and $g_n(t)$ is
linear in $[1/2^{n+2}-1/n2^{n+2},1/2^{n+2}-1/2n2^{n+2}]$ and
$[1/2^n+1/2n2^{n+2},1/2^n+1/n2^{n+2}],$ $n=1,2,....$

Define $e_n=f_{n+1}(a)-f_n(a)$ and $d_1=f_1(a),$ $ d_n=g_n(a),$
$n=2,3,....$ One has that
\beq\label{Tproj-n1}
e_nd_n=d_ne_n=e_n, d_nd_m=d_md_n=0\,\,\,{\rm if}\,\,\, |n-m|\ge 2,
n,m=1,2,....
\eneq

Suppose that $H_1$ and $H_2$ are two Hilbert $A$-modules and
$\phi\in B(H_1, H_2)$ is surjective. Suppose also that there is
$\psi\in B(H, H_2).$ Denote by $H_3$ the closure of $\psi(H).$ Then
$H_3$ is countably generated.

Let $p$ be the open projection of $K(H_1)$ associated with the
Hilbert submodule ${\rm ker} \phi.$
Let  $\phi': H_1/{\rm ker}\phi \to H_2$ be the one-to-one and onto bounded
module map such that
\beq\label{Tproj-1}
\phi'(\pi(x))=\phi(x)\tforal x\in H.
\eneq
Denote by $\phi^{\sim}$ the inverse of $\phi'$ which is also a
bounded module map. There is $x_i\in H_1$ such that
\beq\label{Tproj-2}
\phi(x_i)=\psi(e_i),\,\,\,i=1,2,....
\eneq

Let $\{p_\af\}$ be an approximate identity for $K({\rm ker} \phi).$
By 2.12 of \cite{Lninj}, one may view $K({\rm ker}\phi)\subset
K(H_1).$ Then, by \ref{app}
\beq\label{Tproj-3}
\|\pi(x_i)\|=\inf_\af \|(1-p_\af)x_i\|,\,\,\,i=1,2,....
\eneq
For any $\ep>0.$ Choose $p_n$ so that
\beq\label{Tproj-4}
\|(1-p_n)x_n\|\le \|\pi(x_n)\| +\ep/2^{n+1}=
\|\phi^{\sim}(e_i)\|+\ep/2^{n+1},\,\,\,n=1,2,....
\eneq

Put $\xi_n=(1-p_n)x_nd_n.$ Note that $\phi(\xi_n)=e_nd_n=e_n,$
$n=1,2,....$ For each $n,$ and $b\in A,$ define
\beq\label{Tproj-2+}
T(f_n(a)b)=\sum_{i=1}^n\xi_ib\tforal b\in A.
\eneq
By applying \ref{estimate}
\beq
\|T(\sum_{i=k}^{n+k}e_ib)\|&\le&  2\max_{k\le i\le n+k}\|\xi_i\|\|
(f_{n+k}(a)-f_n(a))b\|\\\label{Tproj-2+1} &\le &
2(\|\phi^{\sim}\|+\sum_{i=1}^n \ep/2^{i+1})\|(f_{n+k}(a)-f_k(a))b\|.
\eneq

Therefore, since $\{f_m(a): m=1,2,...\}$ forms an approximate
identity for $\overline{aAa},$ for any $b\in \overline{aA},$
\beq\label{Tproj-3-1}
\lim_{k\to\infty}\|\sum_{i=k}^{k+n} \xi_ib\|\le
2(\|\phi^{\sim}\|+1)\lim_{k\to\infty}\|(f_{n+k}(a)-f_k(a))b\|=0.
\eneq
 Thus, one defines, for each $b\in B,$
 \beq\label{Tproj-3+1}
 T(b)=\sum_{n=1}^{\infty}\xi_n b.
 \eneq
  By (\ref{Tproj-2+1}),
 \beq\label{Tproj-3+2}
 \|T(b)\|\le 2(\|\phi^{\sim}\|+\ep)\|b\|\tforal b\in
 {\overline{aA}}.
 \eneq

So $T$ is well-define map in $B(H, H_2).$ One verifies that
\beq\label{Tproj-4n}
\phi\circ T(b)&=&\phi\circ T(\sum_{n=1}^{\infty} e_nb)
=\phi(\sum_{n=1}^{\infty} \xi_n b)\\
&=& \phi(\sum_{i=1}^{\infty} (1-p_n)x_nd_nb)
= \sum_{n=1}^{\infty} \phi(x_n)b\\
&=& \sum_{n=1}^{\infty} \psi(e_n)b
=\sum_{n=1}^{\infty}\psi(e_nb)\\
&=&\psi(b).
\eneq

\end{proof}

\begin{cor}\label{Anp}
Let $A$ be a \CA\, and let $x_1, x_2,..,x_n\in A.$ Suppose that
$H_i=\overline{x_iA}$ and $K(H_i)$ admits a sequential central
approximate identity, $i=1,2,...,n.$ Then
$$H=H_1\oplus
H_2\oplus\cdots \oplus H_n$$
 is a projective Hilbert $A$-module (with
bounded module maps as morphisms).
\end{cor}

%

\vspace{0.4in}

\noindent email: hlin@uoregon.edu
\end{document}